\documentclass[10pt]{amsart}
\usepackage{epsfig,url}
\usepackage{mathdots}
\usepackage{hyperref}
\usepackage{xcolor}
\usepackage{marginnote}

\newtheorem{theorem}{Theorem}[section]
\newtheorem{lemma}[theorem]{Lemma}

\theoremstyle{definition}

\theoremstyle{remark}
\newtheorem{remark}[theorem]{Remark}

\begin{document}

\title[Flexibility of the method of brackets]{The evaluation of a definite integral by the method of brackets illustrating its flexibility}

\author[Ivan Gonzalez et al]{Ivan Gonzalez}
\address{Instituto de F\'{i}sica y Astronom\'{i}a, Universidad de Valparaiso, Gran Breta\~{n}a $1111$, Valparaiso, Chile}
\email{ivan.gonzalez@uv.cl}

\author[]{John Lopez Santander}
\address{Department of Mathematics,
Tulane University, New Orleans, LA 70118}
\email{jlopez12@math.tulane.edu}

\author[]{Victor H. Moll}
\address{Department of Mathematics,
Tulane University, New Orleans, LA 70118}
\email{vhm@math.tulane.edu}

\subjclass[2010]{Primary 33}

\date{\today}

\keywords{Integrals, Gradshteyn and Ryzhik, method of brackets}

\maketitle

\newcommand{\ba}{\begin{eqnarray}}
\newcommand{\ea}{\end{eqnarray}}
\newcommand{\ift}{\int_{0}^{\infty}}
\newcommand{\nn}{\nonumber}
\newcommand{\no}{\noindent}
\newcommand{\lf}{\left\lfloor}
\newcommand{\rf}{\right\rfloor}
\newcommand{\realpart}{\mathop{\rm Re}\nolimits}
\newcommand{\imagpart}{\mathop{\rm Im}\nolimits}

\newcommand{\op}[1]{\ensuremath{\operatorname{#1}}}

\newmuskip\pFqmuskip

\newcommand*\pFq[6][8]{%
  \begingroup 
  \pFqmuskip=#1mu\relax
  \mathcode`\,=\string"8000
  \begingroup\lccode`\~=`\,
  \lowercase{\endgroup\let~}\pFqcomma
  {}_{#2}F_{#3}{\left[\genfrac..{0pt}{}{#4}{#5};#6\right]}%
  \endgroup
}
\newcommand{\pFqcomma}{\mskip\pFqmuskip}

\newtheorem{Definition}{\bf Definition}[section]
\newtheorem{Thm}[Definition]{\bf Theorem}
\newtheorem{Example}[Definition]{\bf Example}
\newtheorem{Lem}[Definition]{\bf Lemma}
\newtheorem{Cor}[Definition]{\bf Corollary}
\newtheorem{Prop}[Definition]{\bf Proposition}
\numberwithin{equation}{section}

\section{Introduction}
\label{sec-introduction}
\setcounter{equation}{0}

The goal of the present work is to describe a flexible method of integration, the so-called 
\texttt{method of brackets}, by discussing the evaluation of the identity
\begin{equation}
\label{main-int}
 \int_{0}^{\infty} \int_{0}^{\infty} x^{\alpha-1} y^{\beta-1} \texttt{Ei}(- x^{2} y) K_{0} \left( \frac{x}{y} \right) \, dx dy= 
 - \frac{1}{12} \frac{ \Gamma^{2} \left( \frac{\alpha + \beta}{3} \right) \Gamma^{2} \left(  \frac{\alpha - 2 \beta}{6} \right) }
{4^{(2 \beta - \alpha)/6} \Gamma \left( \frac{\alpha + \beta}{3} + 1 \right)}. 
\end{equation}
\noindent 
Here $\alpha, \, \beta \in \mathbb{R}$, the $\Gamma$-function is the classical eulerian integral and 
the functions appearing in the integrand are the Bessel $K_{0}$ function 
and the exponential integral $\texttt{Ei}$.  The definition of these functions is given below. The complexity of the 
integrand is chosen just to illustrate the power of the method of brackets. We present a variety of ways to 
use this method in the proof of \eqref{main-int}.  

\medskip

The evaluation of definite integrals is one of the basic techniques found in the elementary basic Calculus courses.  One of the central problems is to 
create a class of functions $\mathcal{F}$ and decide if the functions in this class may be integrated. At the elementary level, one starts 
with simple powers $x^{n}, \, n \in \mathbb{N}$ in the class $\mathcal{F}$ and requires some algebraic properties on this class. For instance, it is 
often assumed that the class $\mathcal{F}$ should be closed under elementary operations; that is, $\mathcal{F}$ should be closed under addition and products. 
It follows that $\mathcal{F}$ contains all polynomials in $x$. 

The evaluation of 
\begin{equation}
I(n;a,b) = \int_{a}^{b}  x^{n} \, dx  \quad \textnormal{for} \,\, n \in \mathbb{N}, \,\, \textnormal{ and } \, a, \, b \in \mathbb{R}, 
\end{equation}
\noindent
is elementary  since $\mathcal{F}$  contains the primitive of the integrand: 
\begin{equation}
\label{int-1}
\frac{d}{dx} \left( \frac{1}{n+1} x^{n+1} \right) = x^{n}, \quad \textnormal{when } \, n \neq -1.
\end{equation}
\noindent
The formula 
\begin{equation}
I(n;a,b) = \frac{1}{n+1} \left( b^{n+1} - a^{n+1} \right)
\label{int-2}
\end{equation}
is now immediate.  The linearity of integration now shows that the integration of any polynomial can be completed within the 
class $\mathcal{F}$.  

The extension of \eqref{int-2} to a wider range of parameters $n$ is now a question of elementary analysis. Once $x^{\alpha}$ has been 
defined for $\alpha \in \mathbb{R}$, then \eqref{int-2} is valid with $n$ replaced by $\alpha$. Naturally the exception $\alpha = -1$ remains. 
The integration of this final case  requires the introduction of a new function:
\begin{equation}
\ln x := \int_{1}^{x} t^{-1} \, dt,
\end{equation}
\noindent
the classical (natural) logarithm.  This function is now added to the class $\mathcal{F}$ and the process maybe continued. 

Throughout history the evaluated integrals have been collected in tables such as those created by D.~Bierens de Haan \cite{bierens-1867a} and 
extended in the current table by I.~S.~Gradshteyn and I.~M.~Ryzhik \cite{gradshteyn-2015a}. An effort has been made by the community to 
make the entries in these tables be free of errors. This is a continuing task. The techniques developed to prove these evaluations have generated 
a large number of papers and books.  These include olders volumes \cite{bierens-1862a}, some elementary treatises
\cite{moll-2012a,zwillinger-1992a} and \cite{bronstein-1997a} at a more advanced level. 

During the last few years the authors have developed a method of integration, named \texttt{the method of brackets}. It 
consists on a small number of rules described in Section \ref{1example-sec2}, and is based on the expansion of the integrand in a series of the form 
\begin{equation}
f(x) = \sum_{n=0}^{\infty} a_{n}x^{\alpha n + \beta-1} \quad \textnormal{with} \,\, a_{n}, \, \alpha, \, \beta \in \mathbb{R},
\end{equation}
(the extra $-1$ in the exponent simplifies the appearance of some formulas below).  

\begin{remark}
The original ideas for this method of integration came from the analysis of the so-called \texttt{method of negative dimensional integration}, \cite{anastasiou-2000a,halliday-1987a, suzuki-1998b}. 
A detailed discussion of the power  of this method is presented in \cite{vandeusen-2021a}, with a 
manuscript in preparation.
\end{remark}

The goal of this 
note is to describe the flexibility of the method of brackets  by presenting several ways to evaluate the integral
\eqref{main-int}. 
 It is remarkable that the method of brackets evaluates 
this integral, in view of the fact that  both functions in the integrand have logarithmic singularities at the origin. 

\smallskip

The functions appearing in the integrand of \eqref{main-int} are now defined. 

\smallskip

\noindent
\underline{\texttt{The Bessel function $K_{0}(x)$.}} The traditional manner to introduce
 Bessel functions is as solutions of the differential equation
\begin{equation}
x^{2} \frac{d^{2}y}{dx^{2}} + x \frac{dy}{dx} + (x^{2} - \nu^{2}) y = 0.
\end{equation}
\noindent
With appropriate initial conditions, this gives  the function $J_{\nu}$, with  power series expansion
\begin{equation}
\label{bessel-1}
J_{\nu}(x) =  \sum_{k=0}^{\infty} \frac{(-1)^{k}}{k! \Gamma(\nu+k+1)} \left( \frac{x}{2} \right)^{\nu+ 2k}.
\end{equation}
\noindent
A second (linearly independent)  solution is given by the function 
\begin{equation}
Y_{\nu}(x) = \frac{J_{\nu}(x) \cos(\pi \nu) - J_{-\nu}(x)}{\sin( \pi \nu)}, \quad \nu \not\in \mathbb{Z}. 
\end{equation}
\noindent
(The case when $\nu$ is an integer is treated by a limiting procedure).  References for 
Bessel functions include \cite{andrews-1999a,whittaker-2020a} and the encyclopedic  treatise \cite{watson-1966a}. 

The modified Bessel function $K_{\nu}(x)$ is another family of special functions,  closely related to $J_{\nu}$ and $Y_{\nu}$. 
This function has the expansion
\begin{equation}
K_{0}(x) = - \ln \left( \frac{x}{2} \right) \left( \sum_{k=0}^{\infty} \frac{1}{(k!)^{2}} \left( \frac{x}{2} \right)^{2k} \right) + 
\sum_{k=0}^{\infty} \psi(k+1) \frac{1}{(k!)^{2}} \left( \frac{x}{2} \right)^{2k},
\end{equation}
showing a logarithmic singularity at the origin. This appears in \cite{gradshteyn-2015a}, entry $8.447.3$ as
well as \cite{watson-1966a}, formula WA 95(14).


\smallskip

\noindent
\underline{\texttt{The exponential integral $\texttt{Ei}(x)$.}}
This function  is defined by
\begin{equation}
\text{Ei}(x) = \int_{-\infty}^{x} \frac{e^{t}}{t} \, dt
\end{equation}
\noindent
for  $x<0$. In the case $x>0$, one defines it using the Cauchy principal value
\begin{equation}
\text{Ei}(x) = - \lim\limits_{\epsilon \to 0^{+}}
\left[ \int_{-x}^{-\epsilon} \frac{e^{-t}}{t} \, dt +
\int_{\epsilon}^{\infty} \frac{e^{-t}}{t} \, dt \right],
\end{equation}
\noindent
appearing as entry $\mathbf{3.351.6}$ in \cite{gradshteyn-2015a}.  The series expansion 
\begin{equation}
\text{Ei}(-x) =\gamma + \ln x + \sum_{k=1}^{\infty} \frac{(-x)^{k}}{kk!},
\end{equation}
\noindent
with $\gamma$ the Euler's constant, shows that $\texttt{Ei}$ has a logarithmic singularity at zero.


\section{The method of brackets}
\label{1example-sec2}

This section presents the main rules for the method of brackets, a procedure for the evaluation of definite integrals over the half line $[0,\infty)$. The application of the method consists of a small number of rules, deduced in heuristic form, some of which are placed on solid ground \cite{amdeberhan-2012b}.

Consider the integral
\begin{equation}
    I(f)=\int_0^\infty f(x)\, dx,
\end{equation}
and define the formal symbol 
\begin{equation}
\label{bracket-symbol1}
    \langle a \rangle =\int_0^\infty x^{a-1}\, dx,
\end{equation}
called the \textit{bracket} of $a$. Observe that  the integral on \eqref{bracket-symbol1} is divergent for any 
choice of the parameter $a$. To simplify some computations while operating with brackets, introduce the symbol
\begin{equation}
    \phi_n=\frac{(-1)^n}{\Gamma(n+1)}, 
\end{equation}
called the \texttt{indicator} of $n$, and let $\phi_{i_1,i_2,\ldots, i_r}$, denote
 the product $\phi_{i_1}\phi_{i_2}\cdots \phi_{i_r}$.
 
 \smallskip 

Using this notation, the basic rules for the production and evaluation of a so-called \texttt{bracket series} associated to $I(f)$ are described below (see \cite{gonzalez-2010a,gonzalez-2010b} for details).

\medskip

\textbf{Rule $P_1$.} Assume $f$ has the expansion 
\begin{equation}\label{Rule p1}
    f(x)=\sum_{n=0}^{\infty}{\phi_n a_n x^{\alpha n + \beta - 1}}.
\end{equation}
Then $I(f)$ is assigned the \textit{bracket series} 
\begin{equation}
    I(f)=\sum_{n=0}^{\infty}{\phi_n a_n \langle \alpha n + \beta}\rangle.
\end{equation}

\textbf{Rule $P_2$.} For $\alpha\in \mathbb{R}$, the multinomial power $(a_1+a_2+\cdots+a_r)^{\alpha}$ is assigned the $r$-dimensional bracket series
\begin{equation}\label{Rule P2}
    \sum_{n_1\ge 0} \sum_{n_2\ge 0}\cdots \sum_{n_r\ge 0} \phi _{n_1,n_2,\ldots, n_r} a_1^{n_1}\cdots a_r^{n_r}\frac{\langle n_1+\cdots+n_r-\alpha\rangle}{\Gamma(-\alpha)}.
\end{equation}

\textbf{Definition.} 
Each representation of an integral by a bracket series has associated a \textit{complexity index of the representation} via
\begin{equation} \label{Rule P3}
    \text{complexity index= number of sums -number of brackets}.
\end{equation}
It is important to note that the complexity index is attached to each specific representation of the integral and not just to integral itself. The level of difficulty in the analysis of the resulting bracket series increases with the 
complexity index, however it has been shown that the value of integral is independent of the bracket series
 representation.

\medskip

\textbf{Rule $E_1$.}
Let $\alpha, \, \beta \in \mathbb{R}$. The one-dimensional bracket series is assigned the value
\begin{equation} \label{Rule E1}
    \sum_{n=0}^{\infty}{\phi_n a_n \langle \alpha n + \beta}\rangle=\frac{1}{|\alpha|} f(n^*)\Gamma(-n),
\end{equation}
where $n^*$ is obtained from the vanishing of the bracket; that is $n^*$ solves $\alpha n+ \beta=0$. This is precisely the Ramanujan's Master Theorem (see \cite{amdeberhan-2012b}).

\smallskip

The following rule provides a value for multi-dimensional bracket series of index $0$;  that is, the number of sums is equal to the number of brackets.

\medskip

\textbf{Rule $E_2$.}
Let  $a_{ij} \in \mathbb{R}$. Assuming the matrix $A = (a_{ij})$ is non-singular, then the assignment is
\begin{equation}\label{Rule E2}
\begin{split}
&\sum_{n_1\ge 0} \sum_{n_2\ge 0}\cdots \sum_{n_r\ge 0} \phi_{n_1 \cdots n_r} f(n_1,\ldots,n_r)\langle a_{11}n_1 + \cdots + a_{1r}n_r + c_1 \rangle \\
&\hspace{0.5cm}\cdots \langle a_{r1}n_1 + \cdots + a_{rr}n_r + c_r \rangle = \frac{1}{|\det(A)|}f({n_1}^*,\ldots,{n_r}^*)\Gamma({-n_1}^*)\cdots\Gamma({-n_r}^*),
\end{split}
\end{equation} 
where $\{n_i ^* \}$ is the (unique) solution of the linear system obtained from the vanishing of the brackets. There is no assignment if A is singular.

\medskip

\textbf{Rule $E_3$.} The value of a multi-dimensional bracket series of non-negative 
 index is obtained by computing all the contributions of maximal rank using 
  Rule $E_2$. These contributions to the integral appear as series in the free parameters.  There is no assignment to a bracket series of negative index. 
 
 \smallskip
 
 The next result will be used in later sections. 


\begin{lemma} \label{factoring constant from bracket}
  For any $\alpha, \beta \in \mathbb{R}$, with $\alpha\neq 0$, the bracket satisfies 
\begin{equation}
    \langle \alpha \gamma+\beta\rangle=\frac{1}{|\alpha|}\left \langle\gamma+\frac{\beta}{\alpha}\right \rangle.
\end{equation}  
\end{lemma}

\section{Integral representations}
\label{1example-sec3}

The brackets series of an integral may be obtained if the components of the integrand have proper integral 
representations; that is, representations in terms of functions with series expansions. In the case 
considered here the analysis begins with 
\begin{equation}
\label{k0-intdef}
K_{0}(\xi)  =  \int_{0}^{\infty} \frac{\cos( \xi t )}{(t^{2}+1)^{1/2}} \, dt
\end{equation}
\noindent
and 
\begin{equation}
\texttt{Ei}(-\xi) = - \int_{\xi}^{\infty} \frac{\exp(-t)}{t} \, dt,
\label{ei-intdef}
\end{equation}
\noindent
appearing in \cite{abramowitz-1972a} and  \cite{gradshteyn-2015a}, respectively. 

The classical Taylor series 
\begin{equation}
\cos(\xi t) = \sum_{n} \phi_{n} \frac{\Gamma( \tfrac{1}{2} ) \xi^{2n}}{\Gamma\left( \tfrac{1}{2} + n \right) 4^{n}} t^{2n} 
\end{equation}
\noindent
and the expansion 
\begin{equation}
\frac{1}{(t^{2}+1)^{1/2} } = \sum_{m} \sum_{\ell} \phi_{m,\ell} \frac{\langle \tfrac{1}{2} + m + \ell \rangle }{\Gamma \left( \tfrac{1}{2} \right)} t^{2m},
\end{equation}
\noindent
obtained from \eqref{Rule P2}, it follows from \eqref{k0-intdef} that
\begin{eqnarray}
\label{rep-K0}
K_{0}(\xi) & = & \sum_{n} \sum_{m} \sum_{\ell} \phi_{n,m,\ell} \frac{\langle \tfrac{1}{2} + m + \ell \rangle}{4^{n} \Gamma \left( \tfrac{1}{2} + n \right)} 
\xi^{2n} \int_{0}^{\infty} t^{2n+2m} \, dt \\
& = & \sum_{n} \sum_{m} \sum_{\ell} \phi_{n,m,\ell} \frac{\langle \tfrac{1}{2} + m + \ell \rangle \langle 2n+2m+1 \rangle}{4^{n} \Gamma \left( \tfrac{1}{2} + n \right)} 
\xi^{2n}\nonumber \\
& = & \frac{1}{2} 
\sum_{n} \sum_{m} \sum_{\ell} \phi_{n,m,\ell} \frac{\langle \tfrac{1}{2} + m + \ell \rangle \langle n+m+\tfrac{1}{2}  \rangle}{4^{n} \Gamma \left( \tfrac{1}{2} + n \right)} 
\xi^{2n}, \nonumber
\end{eqnarray}
where the last equality follows from Lemma \ref{factoring constant from bracket}.  This is a bracket series of 
$K_{0}(\xi)$. This is an example of how to assign a bracket series to functions without a power series around zero. 

\medskip

In order to obtain a bracket series representation of the exponential integral, start with 
\begin{equation}
\texttt{Ei}(-\xi) = - \int_{\xi}^{\infty} \frac{\exp(-t) }{t} dt = - \int_{0}^{\infty} \frac{\exp[-(z+\xi) ] }{z+\xi } \, dz. 
\end{equation}
\noindent
Expanding the exponential and using the binomial theorem yield
\begin{eqnarray}
\texttt{Ei}(-\xi) & = &  - \sum_{i} \phi_{i} \int_{0}^{\infty} (z + \xi)^{i-1} \, dz \\
& = & - \sum_{i} \sum_{j} \sum_{k} \phi_{i,j,k} \xi^{k} \frac{ \langle 1 - i + j + k \rangle}{\Gamma(1-i)} \int_{0}^{\infty} z^{j} \, dz \nonumber \\
& = & - \sum_{i} \sum_{j} \sum_{k} \phi_{i,j,k} \xi^{k} \frac{ \langle 1 - i + j +k \rangle}{\Gamma(1-i)} \langle j+1 \rangle. \nonumber 
\end{eqnarray}
\noindent
The bracket sum corresponding to the index $j$ can be evaluated. This is because, using the standard procedure, the 
equations appearing for   any combination of free indices 
are satisfied only for  $j=-1$. This leads to the representation
\begin{equation}
\label{rep-Ei}
\texttt{Ei}(- \xi) = - \sum_{i} \sum_{k} \phi_{i,k} \xi^{k} \frac{\langle k - i \rangle }{\Gamma(1-i)}. 
\end{equation}

Let $I(\alpha, \beta)$ denote the integral in \eqref{main-int}. Replacing the representations 
\eqref{rep-K0} and \eqref{rep-Ei} in \eqref{main-int} give 
\begin{multline}
I(\alpha,\beta) = - \frac{1}{2} \sum_{i} \sum_{k} \sum_{n} \sum_{m} \sum_{\ell} \phi_{i,k,n,m,\ell}  \\
\frac{ \langle k - i \rangle  \langle m+ \ell + \tfrac{1}{2} \rangle \langle n + m + \tfrac{1}{2} \rangle }
{ 4^{n} \Gamma(1-i) \Gamma \left( \tfrac{1}{2} + n \right)}
\langle \alpha + 2n + 2k \rangle 
\langle \beta - 2n + k \rangle. 
\end{multline}

The number of sums equals the number of brackets, therefore this is a representation of index zero. 
According to \eqref{Rule E2}, the value of $I(\alpha, \beta)$ is then given by 
\begin{equation}
I(\alpha,\beta) = - \frac{1}{2 | \det A |} \frac{\Gamma(-i^{*} ) \Gamma(-k^{*}) \Gamma(-n^{*}) \Gamma(-m^{*}) \Gamma(-\ell^{*}) }
{4^{n^{*}} \Gamma(1 - i^{*}) \Gamma \left( \tfrac{1}{2} + n^{*} \right)}
\end{equation}
\noindent
where the set $(i^{*}, k^{*}, n^{*}, m^{*}, \ell^{*})$ is the unique solution of 
\begin{equation}
\begin{pmatrix} 
-1 & 1 & 0 & 0 & 0 \\
0 & 0 & 0 & 1 & 1 \\
0 & 0 & 1 & 1 & 0 \\
0 & 2 & 2 & 0 & 0 \\
0 & 1 & -2 & 0 & 0 
\end{pmatrix}
\begin{pmatrix} i \\ k \\ n \\ m \\ \ell \end{pmatrix} = 
\begin{pmatrix} 0 \\ - \tfrac{1}{2} \\ - \tfrac{1}{2} \\ - \alpha \\ - \beta \end{pmatrix}. 
\end{equation}
\noindent
The matrix of coefficients of the system above is denoted by $A$.   The solutions are given by 
\begin{equation*} 
i^{*} = k^{*} = -\tfrac{1}{3}(\alpha + \beta), \,\, n^{*} = \ell^{*} = - \tfrac{1}{6}(\alpha - 2 \beta), \,\, m^{*} = \tfrac{1}{6}(\alpha - 2 \beta) - \tfrac{1}{2},
\end{equation*}
\noindent
and, since $\det(A)=6$, this gives 
\begin{equation}\label{main int result}
I(\alpha, \beta) = - \frac{1}{12} \frac{ \Gamma^{2} \left( \frac{\alpha + \beta}{3} \right) \Gamma^{2} \left(  \frac{\alpha - 2 \beta}{6} \right) }
{4^{(2 \beta - \alpha)/6} \Gamma \left( \frac{\alpha + \beta}{3} + 1 \right)}. 
\end{equation}
\noindent
This is the first proof of \eqref{main-int}.

\begin{remark}
An alternative evaluation of $I(\alpha,\beta)$ is obtained by replacing the integral representations of $K_{0}$ and $\texttt{Ei}$ and write 
\begin{equation*}
I(\alpha,\beta) = - \int_{0}^{\infty}  \int_{0}^{\infty}  \int_{0}^{\infty}  \int_{0}^{\infty} 
x^{\alpha-1} y^{\beta - 1} 
\frac{\exp[-(z+x^{2}y)]}{(z+x^{2}y)} \frac{\cos(xt/y)}{(t^{2}+1)^{1/2}} dt \, dx \, dy \, dz.
\end{equation*}
\noindent
Then evaluate this integral by the method of brackets. The reader is invited to complete the details.
\end{remark}

\section{The use of null and divergent series for the integrand}
\label{1example-nulldiv}
This section presents the evaluation of the integral \eqref{main-int} by using a technique derived from the method of brackets in \cite{gonzalez-2017a}, which allows to assign a (non-classical) series to functions without a regular power series at zero. These new type of  series are classified as divergent and null.  These terms are illustrated next. 

For the functions in the integrand of \eqref{main-int}, the assigned series are (see \cite{gonzalez-2017a} for details):
\begin{equation}
  K_{0}(\xi) =
    \begin{cases}
        \frac{1}{2}\sum\limits_n \phi_n \frac{\Gamma(-n)}{4^n}\xi^{2n}& \text{Divergent representation}\\
        \frac{1}{\xi}\sum\limits_n \phi_n \frac{4^n\Gamma^2\left(n+\frac{1}{2}\right)}{\Gamma(-n)}\xi^{-2n} & \text{Null representation}
    \end{cases}
\end{equation}

\begin{equation}
    \texttt{Ei}(-\xi)=\sum_l \phi_l\frac{\xi^l}{l}\,\,\,\text{Divergent representation.}
\end{equation}

\begin{remark}
The first representation of $K_{0}(\xi)$ is called divergent since all the coefficients are infinite, since they come
from evaluating the gamma function at its poles.  The second representation is called null, since all the 
coefficients vanish. It is clear that these are not classical power series representations. Nevertheless it will be 
seen that these expansion can be used in the method of brackets to evaluate the integral given in  \eqref{main-int}.
\end{remark}

This gives a new approach for the evaluation of the integral \eqref{main-int}. First, considering the divergent versions for $\texttt{Ei}(-x^2y)$ and $K_0(\frac{x}{y})$,
\begin{equation}
\begin{split}
I(\alpha, \beta) &= \int_{0}^{\infty} \int_{0}^{\infty} x^{\alpha-1} y^{\beta-1} \left(\sum_l\phi_l\frac{(x^2y)^l}{l}\right)\left(\frac{1}{2}\sum_n\phi_n\frac{\Gamma(-n)}{4^n}\left(\frac{x}{y}\right)^{2n}\right) \, dx dy\\
&=\frac{1}{2}\sum_l \sum_n \phi_{l,n}\frac{\Gamma(-n)}{4^n l} \int_{0}^{\infty} \int_{0}^{\infty} x^{\alpha+2l+2n-1}y^{\beta+l-2n-1} \, dx dy,
\end{split}
\end{equation}
\noindent
which gives the bracket series 
\begin{equation}
    I(\alpha, \beta)= \frac{1}{2}\sum_l \sum_n \phi_{l,n} \frac{\Gamma(-n)}{4^n l} \langle \alpha +2l+2n\rangle \langle \beta+l-2n\rangle.
\end{equation}
Evaluating this bracket series according to Rule \ref{Rule E2}, it can be shown that this approach gives the same result found in \eqref{main int result}.

On the other hand, using the divergent representation for $\texttt{Ei}(-x^2y)$ and the null representation for $K_0(\frac{x}{y})$, 
\begin{equation}
\begin{split}
I(\alpha, \beta) &= \int_{0}^{\infty}\hspace{-0.2cm}\int_{0}^{\infty} x^{\alpha-1} y^{\beta-1}\hspace{-0.1cm} \left[\sum_l\phi_l\frac{(x^2y)^l}{l}\right]\hspace{-0.1cm}\left[\sum_n\phi_n\frac{4^n\Gamma^2\left(n+\frac{1}{2}\right)}{\Gamma(-n)}\left(\frac{x}{y}\right)^{-2n-1}\right] dx dy\\
&=\sum_l \sum_n \phi_{l,n}\frac{4^n\Gamma^2\left(n+\frac{1}{2}\right)}{l\Gamma(-n)} \int_{0}^{\infty} \int_{0}^{\infty} x^{\alpha+2l-2n-2}y^{\beta+l+2n} \, dx dy,
\end{split}
\end{equation}
\noindent
which leads to the bracket series 
\begin{equation}
    I(\alpha, \beta)= \sum_l \sum_n \phi_{l,n} \frac{4^n\Gamma^{2}\left(n+\frac{1}{2}\right)}{l\Gamma(-n)}  \langle \alpha+2l-2n-1\rangle \langle \beta+l+2n+1\rangle.
\end{equation}
Now it can be shown that the evaluation of this bracket series gives the same result for  \eqref{main int result} 
as before. 
\section{The use of Mellin transforms}
\label{1example-mellin}

This section illustrates how to compute an integral by using the method of brackets, combined with the fact that
 one may obtain a power series representation for a function whose Mellin transform is known \cite{gonzalez-2020a}. 
 

The main idea is explained next. Let $f(\xi)$ be an arbitrary function with an explicit expression for its Mellin transform
\begin{equation}
    M(s)=\int_0^{\infty} \xi^{s-1}f(\xi)\,d\xi.
\end{equation}
To determine a (perhaps non-classical) power series for $f(\xi)$, suppose 
\begin{equation}
    f(\xi)=\sum_n \phi_n F(n)\xi^{an+b},
\end{equation}
where $a$ and $b$ are real parameters. The coefficients $F(n)$ are determined using 
 the method of brackets:
\begin{equation}
\begin{split}
M(s)&=\int_0^{\infty} \xi^{s-1}f(\xi)\,d\xi\\
&= \sum_n \phi_n F(n)\int_0^{\infty} \xi^{s+an+b-1} \,d\xi\\
&= \sum_n \phi_n F(n)\langle s+an+b \rangle\\
&= \frac{1}{|a|}\Gamma(-n)F(n)\Big|_{n=-\frac{s+b}{a}}.
\end{split}
\end{equation}
This equivalent to 
\begin{equation}
    M(s)\Big|_{s=-an-b}= \frac{1}{|a|}\Gamma(-n)F(n),
\end{equation}
\noindent
which provides an expression for the coefficients $F(n)$ in terms of the Mellin transform of $f$ as 
\begin{equation}
    F(n)=|a|\frac{M(-an-b)}{\Gamma(-n)}.
\end{equation}
Therefore, the series assigned to $f(\xi)$ is
\begin{equation}
    f(\xi)=|a|\sum_n \phi_n \frac{M(-an-b)}{\Gamma(-n)} \xi^{an+b},
\end{equation}
where the parameters $a$ and $b$ are, up to now, arbitrary.

These ideas are now used to calculate the integral \eqref{main-int}. Start with the Mellin transforms
\begin{equation}
    \int_0^{\infty} \xi^{s-1}\texttt{Ei}(-\xi)\,d\xi=-\frac{\Gamma(s)}{s},
\end{equation}
and 
\begin{equation}
    \int_0^{\infty} \xi^{s-1}K_0(\xi)\,d\xi=\frac{2^s}{4}\Gamma^{2}\left(\frac{s}{2}\right),
\end{equation}
appearing as  entries 6.223 and 6.561.16 in \cite{gradshteyn-2015a}.

The computation of the series assigned to the integrand of \eqref{main-int}, starts by writing
\begin{equation}
   \texttt{Ei}(-\xi) =\sum_l \phi_l F(l)\xi^{al+b},
\end{equation}
\begin{equation}
  K_0(\xi)=  \sum_n \phi_n G(n)\xi^{An+B},
\end{equation}
where $a$, $b$, $A$ and $B$ are real. Using the procedure described above, the coefficients of 
$F(l)$ and $G(n)$ are given by
\begin{equation}
  F(l) =|a|\frac{\Gamma(-al-b)}{(al+b)\Gamma(-l)},
\end{equation}
\begin{equation}
 G(n)=2^{-An-B}\frac{|A|}{4}\frac{\Gamma^{2}\left(-\frac{An+B}{2}\right)}{\Gamma(-n)}.
\end{equation}
Thus, the series representations obtained by this method are
\begin{equation} \label{Ei series param ab}
   \texttt{Ei}(-\xi) =|a|\sum_l \phi_l \frac{\Gamma(-al-b)}{(al+b)\Gamma(-l)} \xi^{al+b},
\end{equation}
and
\begin{equation}\label{K_0 series param AB}
  K_0(\xi)= \frac{|A|}{2^{2+B}} \sum_n \phi_n 2^{-An}  \frac{\Gamma^{2}\left(-\frac{An+B}{2}\right)}{\Gamma(-n)} \xi^{An+B}.
\end{equation}
Finally, this yields 
\begin{equation}
\begin{split}
I(\alpha, \beta) &= \int_{0}^{\infty} \int_{0}^{\infty} x^{\alpha-1} y^{\beta-1} \texttt{Ei}(-x^2y) K_0\left(\frac{x}{y}\right) \, dx dy\\
&=\frac{|A||a|}{2^{2+B}}\sum_l \sum_n \phi_{l,n}\frac{\Gamma(-al-b)\Gamma^{2}\left(-\frac{An+B}{2}\right)}{2^{An}(al+b)\Gamma(-l)\Gamma(-n)} \\
&\qquad\qquad\qquad\quad\times \int_{0}^{\infty} \int_{0}^{\infty} x^{\alpha+2al+2b+An+B-1}y^{\beta+al+b-An-B-1}\, dx dy\\
&=\frac{|A||a|}{2^{2+B}}\sum_l \sum_n \phi_{l,n}\frac{\Gamma(-al-b)\Gamma^{2}\left(-\frac{An+B}{2}\right)}{2^{An}(al+b)\Gamma(-l)\Gamma(-n)} \\
&\qquad\qquad\qquad\quad\times \langle \alpha+2al+2b+An+B\rangle \langle\beta+al+b-An-B\rangle,
\end{split}
\end{equation}
\noindent
and then use Rule \ref{Rule E2} to evaluate the  bracket series and  provide a new proof of \eqref{main int result}. 

 

\section{The use of contour integrals}
\label{1example-contour}
This section presents how the method of brackets extends naturally to Mellin-Barnes integrals \cite{gonzalez-2022a}. This technique is based in the following rule for the evaluation of a \textit{bracket integral} over a complex contour.

\medskip

\textbf{Rule $E_4$.}
Assume $F$ is a function defined on $\mathbb{C}$, Then 
\begin{equation}\label{rule eval complex}
    \int_{-i\infty}^{i\infty} F(s)\langle \alpha +\beta s\rangle\, ds =\frac{2\pi i}{|\beta|}F\left(-\frac{\alpha}{\beta}\right).
\end{equation}
This rule has a multi-dimensional version described below. 

\textbf{Rule $E_5$.}
Let $A = (a_{ij})$ be a non-singular matrix, the following expression
\begin{equation}
\begin{split}
    I = & \left(\frac{1}{2\pi i}\right)^N \int_{-i\infty}^{i\infty}\cdots \int_{-i\infty}^{i\infty} F(s_1,\ldots,s_N) \\
    & \times \langle a_{11}s_1+\cdots a_{1N}s_N + c_1\rangle \cdots\langle a_{N1}s_1+\cdots+ a_{NN}s_N + c_N\rangle  \prod_{k=1}^N ds_k
\end{split}
\end{equation}
is evaluated as the value
\begin{equation}\label{rule eval mult complex}
    I=\frac{1}{|\det(A)|}F(s_1^*,\ldots,s_N^*),
\end{equation}
where $\{s_k^*\}$ is the (unique) solution of the linear system obtained from the vanishing of the brackets.

\medskip

The rules presented above are now used to evaluate the integral of interest \eqref{main-int}. Recall that
\begin{equation}
    \int_0^{\infty} \xi^{s-1}\texttt{Ei}(-\xi)\,d\xi=-\frac{\Gamma(s)}{s}
\end{equation}
and 
\begin{equation}
    \int_0^{\infty} \xi^{s-1}K_0(\xi)\,d\xi=\frac{2^s}{4}\Gamma^{2}\left(\frac{s}{2}\right).
\end{equation}
This yields the Mellin-Barnes representations:
\begin{equation}
   \texttt{Ei}(-\xi) =-\frac{1}{2\pi i}\int_{-i\infty}^{i\infty} \xi^{-s}\frac{\Gamma(s)}{s}\, ds
\end{equation}
\begin{equation}
  K_0(\xi)= \frac{1}{8\pi i}\int_{-i\infty}^{i\infty} \xi^{-z} 2^z\Gamma^{2}\left(\frac{z}{2}\right)\, dz.
\end{equation}
Replacing this in \eqref{main-int} gives
\begin{equation}
\begin{split}
I(\alpha, \beta) &= \int_{0}^{\infty} \int_{0}^{\infty} x^{\alpha-1} y^{\beta-1} \texttt{Ei}(-x^2y) K_0\left(\frac{x}{y}\right) \, dx dy\\
&= \int_{0}^{\infty} \int_{0}^{\infty} x^{\alpha-1} y^{\beta-1} \left[-\frac{1}{2\pi i}\int_{-i\infty}^{i\infty} (x^2y)^{-s}\frac{\Gamma(s)}{s}\, ds\right] \\
&\hspace{4.5cm}\times \left[\frac{1}{8\pi i}\int_{-i\infty}^{i\infty} \left(\frac{x}{y}\right)^{-z} 2^z\Gamma^{2}\left(\frac{z}{2}\right)\, dz\right]\, dx dy\\
&= \frac{-1}{4(2\pi i)^2} \int_{-i\infty}^{i\infty} \int_{-i\infty}^{i\infty}2^z\frac{\Gamma(s)}{s}\Gamma^{2}\left(\frac{z}{2}\right) \left[\int_{0}^{\infty} \int_{0}^{\infty} x^{\alpha-2s-z-1}y^{\beta-s+z-1}dxdy\right],
\end{split}
\end{equation}
\noindent
from which one obtains  the bracket integral
\begin{equation}
   I(\alpha, \beta)=  -\frac{1}{4}\frac{1}{(2\pi i)^2} \int_{-i\infty}^{i\infty} \int_{-i\infty}^{i\infty}2^z\frac{\Gamma(s)}{s}\Gamma^{2}\left(\frac{z}{2}\right) \langle\alpha-2s-z\rangle \langle\beta-s+z\rangle\,dsdz.
\end{equation}
This is evaluated  by the repeated use of rule \eqref{rule eval complex} on each variable, or via rule \eqref{rule eval mult complex} with the result
\begin{equation} 
   I(\alpha, \beta) =-\frac{1}{4|-3|}2^z\frac{\Gamma(s)}{s}\Gamma^{2}\left(\frac{z}{2}\right) \Big|_{s^*,\,z^*},
\end{equation}
where $s*=\frac{1}{3}\alpha+\frac{1}{3}\beta$ and $z^*=\frac{1}{3}\alpha-\frac{2}{3}\beta$. Simplifying this expression
 leads to another proof of  \eqref{main int result}.

\section{A mixed method}
\label{1example-mixed}

This section shows that the method of brackets can be applied by combining the techniques from previous sections.
Three different approaches will be analyzed.  The reader is encouraged to explore other possible combinations.

\subsection{Divergent series and Mellin-Barnes, version 1}
Combine the divergent series representation for $\texttt{Ei}(-\xi)$ with the Mellin-Barnes representation for $K_0(\xi)$. Recall these representations are:
\begin{equation}
    \texttt{Ei}(-\xi)=\sum_l \phi_l\frac{\xi^l}{l},
\end{equation}
\begin{equation}
  K_0(\xi)= \frac{1}{8\pi i}\int_{-i\infty}^{i\infty} \xi^{-z} 2^z\Gamma^{2}\left(\frac{z}{2}\right)\, dz.
\end{equation}
Then
\begin{equation}
\begin{split}
I(\alpha, \beta) &= \int_{0}^{\infty} \int_{0}^{\infty} x^{\alpha-1} y^{\beta-1} \texttt{Ei}(-x^2y) K_0\left(\frac{x}{y}\right) \, dx dy\\
&= \int_{0}^{\infty} \int_{0}^{\infty} x^{\alpha-1} y^{\beta-1} \left[\sum_l \phi_l\frac{\xi^l}{l}\right] \left[\frac{1}{8\pi i}\int_{-i\infty}^{i\infty} \left(\frac{x}{y}\right)^{-z} 2^z\Gamma^{2}\left(\frac{z}{2}\right)\, dz\right]\, dx dy\\
&= \frac{1}{8\pi i}\int_{-i\infty}^{i\infty} \sum_l \phi_l\frac{2^z\Gamma^{2}\left(\frac{z}{2}\right)}{l} \left[\int_{0}^{\infty} \int_{0}^{\infty} x^{\alpha+2l-z-1}y^{\beta+l+z-1}\,dxdy\right]\,dz\\
&= \frac{1}{8\pi i}\int_{-i\infty}^{i\infty} \sum_l \phi_l\frac{2^z\Gamma^{2}\left(\frac{z}{2}\right)}{l} \langle\alpha+2l-z\rangle \langle\beta+l+z\rangle.
\end{split}
\end{equation}
Now use the rules of the brackets to evaluate the sum and the integral. First use rule \eqref{rule eval complex} to 
produce
\begin{equation}
\begin{split}
I(\alpha, \beta) &= \frac{1}{4}\sum_l \phi_l\frac{2^z\Gamma^{2}\left(\frac{z}{2}\right)}{l} \langle\beta+l+z\rangle \Big |_{z=\alpha+2l}\\
&= \frac{1}{4}\sum_l \phi_l\frac{2^{\alpha+2l}\Gamma^{2}\left(\frac{\alpha+2l}{2}\right)}{l} \langle\beta+\alpha+3l\rangle \\
&= \frac{1}{12}\sum_l \phi_l\frac{2^{\alpha+2l}\Gamma^{2}\left(\frac{\alpha+2l}{2}\right)}{l} \left\langle\frac{\beta+\alpha}{3}+l\right\rangle,
\end{split}
\end{equation}
using  Lemma \ref{factoring constant from bracket} in the last equality.

Then, evaluate  the  bracket series with Rule \eqref{Rule E1} to obtain
\begin{equation}
\begin{split}
I(\alpha, \beta) &= \frac{1}{12} \frac{2^{\alpha+2l}\Gamma^{2}\left(\frac{\alpha+2l}{2}\right)\Gamma(-l)}{l}  \Big |_{l=-\frac{\beta+\alpha}{3}}\\
&= \frac{1}{12} \frac{2^{\alpha+2\left(-\frac{\beta+\alpha}{3}\right)}\Gamma^{2}\left(\frac{\alpha+2\left(-\frac{\beta+\alpha}{3}\right)}{2}\right)\Gamma\left(\frac{\beta+\alpha}{3}\right)}{\left(-\frac{\beta+\alpha}{3}\right)}\\
&= - \frac{1}{12} \frac{ \Gamma^{2} \left( \frac{\alpha + \beta}{3} \right) \Gamma^{2} \left(  \frac{\alpha - 2 \beta}{6} \right) }
{4^{(2 \beta - \alpha)/6} \Gamma \left( \frac{\alpha + \beta}{3} + 1 \right)},
\end{split}
\end{equation}
\noindent
which simplifies to the usual result in \eqref{main int result}

\subsection{Divergent series and Mellin-Barnes, version 2}
Next use the divergent series for $\texttt{Ei}(-\xi)$ and a bracket series for $K_0(\xi)$ obtained from its corresponding Mellin-Barnes representation. Start with the bracket series for the gamma function,
\begin{equation*}
    \Gamma(\xi)=\sum_n\phi_n\langle\xi+n\rangle.
\end{equation*}
Then, the Mellin-Barnes representation for $K_0(\xi)$ gives
\begin{equation}
\begin{split}
K_0(\xi) &= \frac{1}{8\pi i}\int_{-i\infty}^{i\infty} \xi^{-z} 2^z\Gamma^{2}\left(\frac{z}{2}\right)\, dz\\
     &=\frac{1}{4\pi i}\int_{-i\infty}^{i\infty} \xi^{-2t} 4^t\Gamma^{2}(t)\, dt\\
     &= \frac{1}{4\pi i}\int_{-i\infty}^{i\infty} \xi^{-2t} 4^t\left(\sum_n\phi_n\langle t+n\rangle\right)\left(\sum_m\phi_m\langle t+m\rangle\right)\, dt\\
    &= \frac{1}{4\pi i}\sum_{n,m} \phi_{n,m} \int_{-i\infty}^{i\infty} \xi^{-2t} 4^t \langle t+n\rangle \langle t+m\rangle\, dt.
\end{split}
\end{equation}

Using  \eqref{rule eval complex} on the bracket $\langle t+m\rangle$  eliminates the integral to obtain
\begin{equation}\label{K_0 bracket rep 2parameters}
    K_0(\xi)=\frac{1}{2} \sum_{m,n} \phi_{n,m} \xi^{2m} 4^{-m} \langle n-m\rangle,
\end{equation}
which is another bracket series for $K_0(\xi)$. Then the integral of interest is given by
\begin{equation}
\begin{split}
I(\alpha, \beta) &= \int_{0}^{\infty} \int_{0}^{\infty} x^{\alpha-1} y^{\beta-1} \texttt{Ei}(-x^2y) K_0\left(\frac{x}{y}\right) \, dx dy\\
&=\int_{0}^{\infty} \int_{0}^{\infty} x^{\alpha-1} y^{\beta-1} \left(\sum_l \phi_l\frac{(x^2y)^l}{l}\right)\left(\frac{1}{2} \sum_{m,n} \phi_{n,m} \left(\frac{x}{y}\right)^{2m}4^{-m} \langle n-m\rangle\right)\, dx dy\\
&=\sum_{l,n,m}\phi_{l,n,m}\frac{4^{-m}\langle n-m\rangle}{l} \int_{0}^{\infty} \int_{0}^{\infty}x^{\alpha+2l+2m-1}y^{\beta+l-2m-1}\,dxdy,
\end{split}
\end{equation}
\noindent
which yields the three-dimensional bracket series
\begin{equation}
    I(\alpha, \beta)=\sum_{l,n,m}\phi_{l,n,m}\frac{\langle n-m\rangle \langle \alpha+2l+2m\rangle \langle \beta+l-2m \rangle}{4^m l}.
\end{equation}
Rule \eqref{Rule E2} evaluates this bracket series and gives the usual result for \eqref{main int result}.

\begin{remark}
    The bracket series representation for $K_0(\xi)$ obtained in \eqref{K_0 bracket rep 2parameters}, can be evaluated using rule \eqref{Rule E1} to eliminate one of the two sums. Independently of the choice of the parameter to be eliminated, it gives
\begin{equation}
    K_0(\xi)=\frac{1}{2} \sum_{k} \phi_{k} \frac{\Gamma(-k)}{4^k}\xi^{2k},
\end{equation} 
the previously known divergent series representation for $K_0(\xi)$. This is obtained now by  a completely different 
approach to the one used in \cite{gonzalez-2017a},  illustrating  the flexibility of the method of brackets. 
\end{remark}

\section{Conclusions}
\label{1example-conclusions}

The method of brackets is an algorithm designed to evaluate definite integrals over a half-line. The method is 
based on a small number of operational rules. Recently, \cite{gonzalez-2017a,gonzalez-2020a}, this 
 method has been extended to deal with functions that lack of a power series representations. In the present work 
 we illustrate these new techniques by evaluating an integral where the integrand is a combination of 
 a Bessel function and the  exponential integral function.




\end{document}